\documentclass[12pt]{amsart}

\usepackage{graphicx}
\usepackage{amsmath}

\usepackage{tikz}

\def\frak{\mathfrak}

\def\S{\mathbb{S}}

\def\P{\mathbb{P}}
\def\R{\mathbb{R}}

\def\cD{\mathcal{D}}

\def\al{\alpha}

\def\ka{\kappa}
\def\la{\lambda}

\def\om{\omega}

\def\Si{\Sigma}

\def\Om{\Omega}

\newcommand{\der}{{\rm d}}

\numberwithin{equation}{section}

\newtheorem{theorem}{Theorem}[section]

\newtheorem{proposition}[theorem]{Proposition}

\newtheorem{definition}[theorem]{Definition}
\newtheorem{example}{Example}[section]
\theoremstyle{remark}

\theoremstyle{remark}

\usepackage{amssymb,stmaryrd}
\usepackage{amscd}

\author{Matthew Randall}
\address{School of Mathematics and Statistics\\
Nanjing University of Information Science and Technology\\
219 Ningliu Road\\
Nanjing, China}
\email{100093@nuist.edu.cn}

\title{Local equivalence of some maximally symmetric $(2,3,5)$-distributions III: An-Nurowski surface rolling on a plane}

\begin{document}

\begin{abstract}
We consider the maximally symmetric $(2,3,5)$-distribution given by the An-Nurowski circle twistor bundle over the product of an An-Nurowski surface and the plane. This circle twistor distribution encodes the configuration space of an An-Nurowski surface rolling without slipping or twisting on the plane. We calculate the vector fields associated to this maximally symmetric $(2,3,5)$-distribution that define a split $\frak{g}_2$ Lie algebra and by projection we obtain an action of $SL(3,\R)$ on the submanifold $\R^2\times \S^1$ (the configuration space without a surface). 
\end{abstract}

\maketitle

\pagestyle{myheadings}
\markboth{Randall}{Local equivalence of some maximally symmetric $(2,3,5)$-distributions III}

\section{Introduction}

In the work of \cite{AN14}, it was discovered that in addition to the configuration space of two spheres or two hyperboloids with curvature ratios $1:9$ or $9:1$ rolling without slipping or twisting over each other, there is also the configuration space of an additional surface $\Si^2$ equipped with the Riemannian metric
\begin{align}\label{anm1}
ds^2=\left(\frac{\kappa c^2-h^2}{2c^3}\right)^2\der h^2+h^2 \der q^2
\end{align}
that has maximal $G_2$ symmetry when rolling without slipping or twising on the plane, where $\kappa$ and $c$ are arbitrary constants with $c$ non-zero. We call a Riemmanian surface with such a metric an An-Nurowski surface. It comes in a family parametrised by two constants and modulo homotheties generated by $c$, three inequivalent classes. This is the main Theorem 3 in \cite{AN14}. In this paper we calculate the split $\frak{g}_2$ Lie algebra of vector fields associated to the $(2,3,5)$-distribution given by the configuration space of the An-Nurowski surface $\Si^2$ rolling without slipping or twisting on $\R^2$. A full list of vector fields is given for the surface with $\ka=0$ and $c=1$, corresponding to the metric $g_{1o}$ in \cite{AN14}. As a corollary, by projecting the vector fields that form $\frak{sl}_{3}\subset \frak{g}_2$ we obtain an action of $SL(3,\R)$ on $\R^2 \times \S^1$, the submanifold of the configuration space without the surface. 

The primary goal of the paper is to calculate the $\frak{g}_2$ Lie algebra of vector fields for this specific rolling example. Most if not all of the results here have already been obtained and are implicit in \cite{AN14}, we just work out a few details. The case of hyperboloid rolling has been dealt with in \cite{r17} and spheres rolling in \cite{r21b}.
First, we establish the metric representative of the conformal structure on this configuration space. 

\section{Nurowski's metric for the configruation space of An-Nurowski surface rolling on a plane} 
Consider the coframing
\begin{align*}
t_1&=\der \phi+\der x,\quad t_2=\der\phi-\der x,\quad t_3=\der \theta+h(x) \der q,\quad t_4=\der \theta-h(x) \der q
\end{align*}
on the product of the plane $\R^2$ with local coordinates $(\phi, \theta)$ 
and standard metric 
\begin{align}\label{plane}
ds_1^2=\der \theta^2+\der \phi^2
 \end{align}
and a surface of revolution $\Si^2$ with coordinates $(x,q)$ and the `warped product' metric
\begin{align}\label{surfrev}
ds_2^2=\der x^2+h(x)^2\der q^2.
 \end{align}
To the product of these two surfaces with the split-signature product metric
\[
2 t_1 t_2+2 t_3 t_4
\]
we can form the five dimensional An-Nurowski circle twistor distribution $M$ which parametrises the self-dual null 2-planes of $\R^2 \times \Si^2$ as a circle bundle. This circle twistor distribution is constructed by the horizontal lifting of the self-dual null 2-planes and is naturally equipped with a $(2,3,5)$-distribution as follows. 

The null 2-planes are annihilated by the two 1-forms
\begin{align*}
\cos(\psi)t_3-(\sin(\psi)+1)t_1\hspace{12pt}\mbox{and}\hspace{12pt} \cos(\psi)t_2+(\sin(\psi)+1)t_4. 
\end{align*}
The requirement that the self-dual null 2-planes are parallel gives the additional 1-form
\begin{align*}
d\psi+\frac{h'}{2h}(t_3-t_4). 
\end{align*}

These three 1-forms encode the An-Nurowski circle twistor distribution associated to a surface of revolution (with Killing symmetry) rolling on the plane without slipping and twisting over each other. 
The distribution $\cD$ spanned by the vector fields 
\begin{align*}
V^1=\partial_{\theta}+\cos(\psi)\partial_x+\sin(\psi)(\frac{1}{h}\partial_q-\frac{h_x}{h}\partial_\psi),\\
V^2=\partial_{\phi}-\sin(\psi)\partial_x+\cos(\psi)(\frac{1}{h}\partial_q-\frac{h_x}{h}\partial_\psi)
\end{align*}
is annihilated by the three 1-forms and has rank $2$. The bracket distribution $[\cD, \cD]$ has rank 3 whenever $h_{xx}\neq 0$ while the full tangent space $[\cD,[\cD,\cD]]\cong TM$ has rank 5, hence such a distribution has growth vector $(2,3,5)$ provided that the genericity condition $h_{xx}\neq 0$ holds. 

In Section 5 of \cite{PN05-conf}, it is shown how to associate canonically to any $(2,3,5)$-distribution a conformal class of metrics of split signature $(2,3)$, also known as Nurowski's conformal structure or Nurowski's conformal metrics, such that the rank 2 distribution is isotropic with respect to any metric in the conformal class. The method of equivalence \cite{Cartan1910} (also see Section 5 of \cite{PN05-conf}) applied to $(2,3,5)$-distributions produces the 1-forms $(\theta_1, \theta_2,\theta_3, \theta_4, \theta_5)$ and $(\Om_1, \ldots, \Om_9)$ that together define a rank 14 principal bundle over the 5-manifold $M$. A representative metric in Nurowski's conformal class is given by
\begin{align}\label{metric}
g=2 \theta_1 \theta_5-2\theta_2 \theta_4+\frac{4}{3}\theta_3 \theta_3.
\end{align}
When $g$ has vanishing Weyl tensor, the distribution is called maximally symmetric and has the smallest exceptional Lie group split $G_2$ as its group of local symmetries. 

Let us now determine this representative metric of the conformal structure on the circle twistor distribution for our example of the surface of revolution rolling on the plane, with the genericity condition $h_{xx}\neq 0$. 
Using the notation $\chi$, we have
 \begin{align*}
 \chi_1=\der \theta,\qquad \chi_2=\der \phi,\qquad \chi_3=\der x,\qquad \chi_4=h(x)\der q.
 \end{align*}
 The distribution $\cD$ is annihilated by the 1-forms
 \begin{align*}
 \om_1&=\chi_1-\cos(\psi)\chi_3-\sin(\psi) \chi_4=\frac{1+\sin(\psi)}{2}t_1+\frac{1-\sin(\psi)}{2}t_2+\frac{\cos(\psi)}{2}(t_3-t_4),\\
 \om_2&=\chi_2+\sin(\psi)\chi_3-\cos(\psi)\chi_4=\frac{1+\sin(\psi)}{2}t_4+\frac{1-\sin(\psi)}{2}t_3+\frac{\cos(\psi)}{2}(t_2-t_1),\\
 \om_3&=\der\psi+\frac{h_{x}}{h}\chi_4. 
 \end{align*}
 This can be completed to a coframe by taking
 \begin{align*}
 \om_4=\sin(\psi)\chi_3-\cos(\psi) \chi_4,\quad \om_5=\cos(\psi)\chi_3+\sin(\psi)\chi_4.
 \end{align*}

To obtain the coframing for Nurowski's conformal structure, we take
 \begin{align*}
 \theta_1&=\om_1,\quad \theta_2=\om_2,\quad \theta_3=\left(\frac{h}{h''}\right)^{1/3}\om_3,\\
 \theta_4&=\left(\frac{h''}{h}\right)^{1/3}\bigg[\om_4+\left(\frac{\la}{10 hh''^{2}}-\frac{\mu \cos(\psi)^2}{30hh'^2h''^3}\right)\om_2+\frac{(3hh''^2+\la)\cos(\psi)}{3h'(h'')^2}\om_3\bigg],\\
 \theta_5&=\left(\frac{h''}{h}\right)^{1/3}\bigg[\om_5-\left(\frac{\la}{10 hh''^{2}}-\frac{\mu\sin(\psi)^2}{30hh'^2h''^{3}}\right)\om_1+\frac{\mu\cos(\psi)\sin(\psi)}{15hh'^2h''^3}\om_2-\frac{(3hh''^2+\la)\sin(\psi)}{3h'(h'')^2}\om_3\bigg],
 \end{align*}
where \[
\la=h'''h'h-h''h'^2-3h''^2h \quad\mbox{and}\quad \mu=15\la h''^2h-3 h' h''(\la' h-3\la h')+5 \la^2.
\]
 From this, we obtain
 \begin{align*}
 g&=2\theta_1\theta_5-2\theta_2\theta_4+\frac{4}{3}\theta_3^2\\
 &=\left(\frac{h''}{h}\right)^{1/3}\bigg[2(\om_1\om_5-\om_2\om_4)-\frac{\la}{5hh''^2}(\om_1^2+\om_2^2)+\frac{\mu}{15hh'^2h''^3}(\sin(\psi)\om_1+\cos(\psi)\om_2)^2\\
&\quad{}-\frac{2(3hh''^2+\la)}{3h'(h'')^{2}}(\sin(\psi)\om_1\om_3+\cos(\psi)\om_2\om_3)\bigg]+\frac{4}{3}\left(\frac{h}{h''}\right)^{2/3}\om_3^2,
 \end{align*}
which can be rescaled by the function $(\frac{h''}{h})^{-1/3}$ to give
 \begin{align*}
 \tilde g&=2(\om_1\om_5-\om_2\om_4)-\frac{\la}{5hh''^2}(\om_1^2+\om_2^2)+\frac{\mu}{15hh'^2h''^3}(\sin(\psi)\om_1+\cos(\psi)\om_2)^2\\
&\quad{}-\frac{2(3hh''^2+\la)}{3h'h''^2}(\sin(\psi)\om_1\om_3+\cos(\psi)\om_2\om_3)+\frac{4h}{3h''}\om_3^2.
 \end{align*}
 Since
 \begin{align*}
 \om_1&=\chi_1-\om_5,\quad{} \om_2=\chi_2+\om_4,\quad{}  \om_4^2+\om_5^2=\chi_3^2+\chi_4^2,\\
 \chi_1\om_5&=\frac{1}{2}(\chi_1^2+\om_5^2-(\chi_1-\om_5)^2),\quad{} \chi_2\om_4=\frac{1}{2}((\chi_2+\om_4)^2-\chi_2^2-\om_4^2),
 \end{align*}
 we have 
 \begin{align*}
 2(\om_1\om_5-\om_2\om_4)=\chi_1^2+\chi_2^2-\chi_3^2-\chi_4^2-\om_1^2-\om_2^2
 \end{align*}
and we can therefore express the metric more simply as
\begin{align*}
 \tilde g&=\chi_1^2+\chi_2^2-\chi_3^2-\chi_4^2-\om_1^2-\om_2^2-\frac{\la}{5hh''^2}(\om_1^2+\om_2^2)+\frac{\mu}{15hh'^2h''^3}(\sin(\psi)\om_1+\cos(\psi)\om_2)^2\\
&\quad{}-\frac{2(3hh''^2+\la)}{3h'h''^2}(\sin(\psi)\om_1\om_3+\cos(\psi)\om_2\om_3)+\frac{4h}{3h''}\om_3^2.
 \end{align*}
It was determined in \cite{AN14} that the metric is conformally flat and the $(2,3,5)$-distribution achieves maximal symmetry when $\la(x)=0$, that is to say $h(x)$ satisfies the third-order ordinary differential equation
 \begin{align*}
 h'''h'h-3h''^2h-h''h'^2=0.
 \end{align*}
 \begin{theorem}[Proposition 5 of \cite{AN14}] Nurowski's representative metric for the conformal structure associated to the An-Nurowski circle twistor distribution over the product of two surfaces, one a surface of revolution with metric (\ref{surfrev}) and $h''\neq 0$ and the other of the plane with metric (\ref{plane}), rolling without slipping or twisting over each other is given by
\begin{align*}
 \tilde g&=\chi_1^2+\chi_2^2-\chi_3^2-\chi_4^2-(\frac{\la}{5hh''^2}+1)(\om_1^2+\om_2^2)+\frac{\mu}{15hh'^2h''^3}(\sin(\psi)\om_1+\cos(\psi)\om_2)^2\\
&\quad{}-\frac{2(3hh''^2+\la)}{3h'h''^2}(\sin(\psi)\om_1\om_3+\cos(\psi)\om_2\om_3)+\frac{4h}{3h''}\om_3^2.
 \end{align*}
 The metric is conformally flat and the distribution has maximal split $G_2$ symmetry whenever the differential equation 
  \begin{align}\label{h3}
 h'''h'h-3h''^2h-h''h'^2=0
 \end{align}
 is satisfied.
 \end{theorem}
\begin{example}[Sphere rolling on a plane]
In the case where we have the configuration space of the unit sphere rolling on the plane, taking $h(x)=\sin(x)$ we find $\mu(x)=0$ and $3hh''^2+\la=0$. This gives the conformal structure
\begin{align*}
 \tilde g&=\chi_1^2+\chi_2^2-\chi_3^2-\chi_4^2-\frac{2}{5}(\om_1^2+\om_2^2)-\frac{4}{3}\om_3^2
 \end{align*}
which is not conformally flat. 
\end{example}

The solutions to equation (\ref{h3}) give rise to the aforementioned surfaces of An and Nurowski. 

\begin{definition}
We shall call a Riemmanian surface of revolution equipped with a metric (\ref{surfrev}), where $h(x)$ is a solution to the differential equaion (\ref{h3}), an {\bf An-Nurowski surface}.  
\end{definition}

\subsection{Relating the differential equations}
 The differential equation (\ref{h3}) is equivalent to the differential equation (4.11) of \cite{AN14}
 \[
 \rho_{xxx}\rho_x\rho^2-3\rho_{xx}^2\rho^2+\rho_{xx}\rho_{x}^2\rho+\rho_x^4=0
 \]
 derived by the authors An and Nurowski under the following change of variables. Take $\tilde x=\int \rho \der x$, where $\rho$, $x$ are coordinates in the proof of Theorem 3 in \cite{AN14}. We take $h(\tilde x)=\rho$. It follows that
 \begin{align*}
 \rho_{x}&=h h_{\tilde x},\\
 \rho_{xx}&=hh_{\tilde x}^2+h^2h_{\tilde x \tilde x},\\
 \rho_{xxx}&=h(h_{\tilde x}^3+4h h_{\tilde x}h_{\tilde x\tilde x}+h^2h_{\tilde x\tilde x\tilde x}).
 \end{align*}
 Substituting these into equation (4.11) of \cite{AN14} gives
 \begin{align*}
 h^5(h_{\tilde x\tilde x\tilde x}h_{\tilde x}h-3h_{\tilde x \tilde x}^2h-h_{\tilde x\tilde x}h_{\tilde x}^2)=0,
 \end{align*}
 and we relabel $\tilde x$ as $x$ to get equation (\ref{h3}).
 
\subsection{Solving the differential equation}
 We now consider the solutions to the differential equation (\ref{h3})
  \begin{align*}
 h'''h'h-3h''^2h-h''h'^2=0.
 \end{align*}
 This can be solved in the following manner. The left hand side is the numerator of the derivative of $\frac{h(h')^3}{h''}$, so one integration of
 \[
 \frac{\der }{\der x}\left(\frac{h(h')^3}{h''}\right)=0
 \]
 gives
 \[
 \frac{h(h')^3}{h''}=c^3.
 \]
 This second-order equation admits one more integral, since this equation is the numerator of the derivative of
 \[
 \frac{h^2}{c^2}+\frac{2c}{h_x}.
 \]
 It follows that 
 \[
  \frac{h^2}{c^2}+\frac{2c}{h_x}=\kappa, 
 \]
 where $\kappa$ and $c$ are two constants of integration, with $c$ non-zero. 
 We therefore have
 \[
 h_x=\frac{2c^3}{\kappa c^2-h^2}.
 \]
This gives
 \[
\der h=\frac{2c^3}{\ka c^2-h^2}\der x, 
\]
or equivalently, 
 \[
\frac{\ka c^2-h^2}{2c^3}\der h=\der x, 
\]
from which we obtain
\begin{align*}
x+\al&=\int\frac{\ka c^2-h^2}{2c^3}\der h\\
&=\frac{\ka}{2c}h-\frac{1}{6c^3}h^3,
\end{align*}
and the solutions to (\ref{h3}) are given by roots of the cubic equation
\[
h^3-3\ka c^2h+6c^3(x+\al)=0.
\]
 We now make the change of variables
 \[
\der h=\frac{2c^3}{\ka c^2-h^2}\der x, 
\]
 replacing $\partial_x$ by $\frac{2c^3}{\ka c^2-h^2}\partial_h$. In this new coordinate, the metric (\ref{surfrev}) becomes that of (\ref{anm1}).

In this new coordinate system, the result about the metric on the circle twistor bundle obtained over the product of an An-Nurowski metric and the plane reads
  \begin{theorem}\label{ans1} Nurowski's representative metric for the conformal structure associated to the An-Nurowski circle twistor distribution over the product of two surfaces, one an An-Nurowski surface with metric (\ref{anm1}) and the other of the plane, rolling without slipping or twisting over each other is given by
\begin{align}\label{ganm}
 \tilde g=\chi_1^2+\chi_2^2-\chi_3^2-\chi_4^2-2\frac{(\ka c^2-h^2)h}{2c^3}(\sin(\psi)\om_1\om_3+\cos(\psi)\om_2\om_3)+\frac{(\ka c^2-h^2)^3}{6c^6}\om_3^2-\om_1^2-\om_2^2.
 \end{align}
 The metric is conformally flat and the distribution has maximal split $G_2$ symmetry.
 \end{theorem}
 In this coordinate system $(\theta, \phi, \psi, h, q)$, we shall now write down the vector fields that bracket-generate the Lie algebra of split $\frak{g}_2$. 
 
 \section{Split $\frak{g}_2$ Lie algebra associated to the An-Nurowski circle twistor distribution}
 In this section, we consider the $(2,3,5)$-distribution \`a la An and Nurowski \cite{AN14} in the case of An-Nurowski surfaces rolling without slipping or twisting over a plane (or the plane rolling without slipping or twisting over the surface). 
 The surface metric given by
 \begin{align*}
 ds_1^2=\der \theta^2+ \der \phi^2
 \end{align*}
is the metric on the flat plane and has constant negative Gauss curvature $0$. 
 The An-Nurowski surface with redefined coordinates $(h,q)$ has the metric given by
 \begin{align*}
 ds_2^2=\left(\frac{\ka c^2-h^2}{2c^3}\right)^2\der h^2+h^2\der q^2
 \end{align*}
with negative Gauss curvature $-\frac{8c^6}{(\ka c^2-h^2)^3}$ and it admits a Killing symmetry. 
 
 To relate the maximally symmetric rolling $(2,3,5)$-distribution to the flat Cartan distribution obtained in \cite{Cartan1910}, we take the coordinate transformation given by
 \begin{align*}
c_1&=6c^3\theta-4h(2c^3q+\ka c^2\psi)\sin(\psi)-h(3c^2\ka-h^2)\cos(\psi),\\
c_2&=6c^3\phi-4h(2c^3q+\ka c^2\psi)\cos(\psi)+h(3c^2\ka-h^2)\sin(\psi),\\
c_3&=2 c^3 q+\ka c^2 \psi,\\
c_4&=-\sin(\psi) h,\\
c_5&=\cos(\psi) h.
 \end{align*}

 It follows that the following 1-forms
\begin{align*}
\Theta_1&:=\der c_1-2c_4 \der c_3-4 c_3 \der c_4,\\
\Theta_2&:=\der c_2+2 c_5 \der c_3+4 c_3\der c_5,\\
\Theta_3&:=\der c_3+c_5 \der c_4-c_4 \der c_5,
\end{align*}
 are in the span of the 1-forms $(\om_1,\om_2,\om_3)$ given by
 \begin{align*}
 \om_1&=\der \theta-\cos(\psi)\frac{\ka c^2-h^2}{2c^3}\der h-\sin(\psi) h \der q,\\
 \om_2&=\der \phi+\sin(\psi)\frac{\ka c^2-h^2}{2c^3}\der h-\cos(\psi) h \der q,\\
  \om_3&=\der\psi+\frac{2c^3}{\ka c^2-h^2}\der q.
\end{align*}
Under the above change of coordinates, we have
\begin{align*}
\Theta_1&=6c^3\om_1+h\sin(\psi)(\ka c^2-h^2)\om_3,\\
\Theta_2&=6c^3\om_2+h\cos(\psi)(\ka c^2-h^2)\om_3,\\
\Theta_3&=-(h^2-\kappa c^2) \om_3.
\end{align*}

The distribution is spanned by the vector fields
\begin{align*}
X^1&=\partial_\theta+\frac{\sin(\psi)}{h}\left(\partial_q-\frac{2c^3}{\ka c^2-h^2}\partial_{\psi}\right)+\frac{2c^3}{\ka c^2-h^2}\cos(\psi)\partial_h,\\
X^2&=\partial_\phi+\frac{\cos(\psi)}{h}\left(\partial_q-\frac{2c^3}{\ka c^2-h^2}\partial_{\psi}\right)-\frac{2c^3}{\ka c^2-h^2}\sin(\psi)\partial_h.
 \end{align*} 

The coordinates for the flat metric are obtained by passing to the Engel distribution by taking 
$\frak{r_1}=c_5$, $\frak{r_2}=c_4$, $\frak{r_3}=c_3$, $\frak{r_4}=\frac{1}{2}(c_2+3c_3c_5)$, $\frak{r_5}=\frac{1}{2}(c_1-3c_3c_4)$.
We find
\begin{align*}
\frak{r_1}&=\cos(\psi) h,\\
\frak{r_2}&=-\sin(\psi) h,\\
\frak{r_3}&=2 c^3 q+\ka c^2 \psi,\\
\frak{r_4}&=-\frac{h}{2}(2c^3q+\ka c^2\psi)\cos(\psi)+3c^3\phi+\frac{h}{2}(3\ka c^2-h^2)\sin(\psi),\\
\frak{r_5}&=-\frac{h}{2}(2c^3q+\ka c^2\psi)\sin(\psi)+3c^3\theta-\frac{h}{2}(3\ka c^2-h^2)\cos(\psi).
\end{align*}
It follows that Nurowski's metric for the maximally symmetric distribution of an An-Nurowski surface rolling on the plane given by (\ref{ganm}) is conformally equivalent to the metric
\begin{align*}
\der s^2=2\der \frak{r_1}\der \frak{r_5}+2\der \frak{r_2}\der \frak{r_4}+\der \frak{r_3}^2
\end{align*}
determined by Engel in \cite{Engel} and \cite{Engel2}. 
\begin{proposition}
The coordinate functions $(\frak{r_1},\frak{r_2},\frak{r_3},\frak{r_4},\frak{r_5})$ given above map the Nurowski conformal structure associated to the maximally symmetric rolling distribution of the An-Nurowski surface over the plane to the standard split signature flat metric on $\R^5$ given by Engel in \cite{Engel} and \cite{Engel2}. We have
\[
\der s^2=\frac{6c^6}{\ka c^2-h^2}\tilde g
\]
where $\tilde g$ is given in (\ref{ganm}) of Theorem \ref{ans1}.
\end{proposition}

 We also have that the 1-forms given in \cite{Engel} and \cite{Engel2} satisfy
\begin{align*}
\der\frak{r_3}+\frak{r_1}\der \frak{r_2}-\frak{r_2}\der \frak{r_1}=\Theta_3,\hspace{4pt} \der\frak{r_4}+\frac{1}{2}(\frak{r_3}\der \frak{r_1}-\frak{r_1}\der \frak{r_3})=\frac{1}{2}\Theta_2,\hspace{4pt} 
\der\frak{r_5}+\frac{1}{2}(\frak{r_2}\der \frak{r_3}-\frak{r_3}\der \frak{r_2})=\frac{1}{2}\Theta_1.
\end{align*}

Following \cite{r21b}, we now write down the vector fields
\begin{align*}
Z^1=\partial_{c_3}+2c_5\partial_{c_2}-2c_4\partial_{c_1},\hspace{4pt} Z^2=\partial_{c_4}+4c_3\partial_{c_1}-2c_5\partial_{c_3},\hspace{4pt} Z^3=\partial_{c_5}+2c_4\partial_{c_3}-4c_3\partial_{c_2},
\end{align*}
and define
\begin{align}\label{cdist}
S^1=Z^2+c_5Z^1,\hspace{12pt} S^2=Z^3-c_4Z^1, \hspace{12pt} S^3=-c_1Z^2+c_2Z^3-(c_1c_5+c_2c_4+c_3^2)Z^1.
\end{align}
The vector fields $S^1$, $S^2$ and $S^3$ pairwise bracket-generate the Lie algebra of split ${\frak g}_2$ in the following sense: 
defining the following vector fields by the Lie brackets
\begin{align*}
S^4&=[S^1,S^2],\quad S^5=[S^2,S^3],\quad S^6=[S^3,S^1],\\
L^1&=[S^1,S^4],\quad L^3=[S^2,S^5],\quad L^5=[S^3,S^6],\\
L^2&=[S^2,S^4],\quad L^4=[S^3,S^5],\quad L^6=[S^1,S^6]
\end{align*}
and taking $H_1=[S^2,S^6]$ and $H_2=[S^4,S^3]$,
we require that the set of vector fields
\begin{align*}
\{S^1,S^2,S^3,S^4,S^5,S^6,\frac{1}{4}(H_2-H_1),\frac{\sqrt{3}}{12}(H_1+H_2),L^1,L^2,L^3,L^4,L^5,L^6\}
\end{align*}
form the 14-dimensional Lie algebra of split ${\frak g}_2$. With this choice of the Cartan subalgebra spanned by $\frac{1}{4}(H_2-H_1)$ and $\frac{\sqrt{3}}{12}(H_1+H_2)$, the root diagram is given by the picture below. 
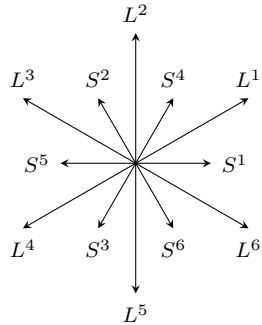
\begin{figure}[h!]
\begin{tikzpicture}
	\draw [stealth-stealth](-1,0) -- (1,0);
\draw (1,0) node[anchor=west] {{\tiny $S^1$}};
\draw (-1,0) node[anchor=east] {{\tiny $S^5$}};
	\draw [stealth-stealth](0,-1.732) -- (0,1.732);
\draw (0,-1.732) node[anchor=north] {{\tiny $L^5$}};
\draw (0,1.732) node[anchor=south] {{\tiny $L^2$}};
\draw [stealth-stealth](-0.5,-0.866) -- (0.5,0.866);
\draw (-0.5,-0.866) node[anchor=north] {{\tiny $S^{3}$}};
\draw (0.5,0.866) node[anchor=south] {{\tiny $S^4$}};
\draw [stealth-stealth](-1.5,-0.866) -- (1.5,0.866);
\draw (-1.5,-0.866) node[anchor=north] {{\tiny $L^4$}};
\draw (1.5,0.866) node[anchor=south] {{\tiny $L^1$}};
\draw [stealth-stealth](1.5,-0.866) -- (-1.5,0.866);
\draw(1.5,-0.866) node[anchor=north] {{\tiny $L^6$}};
\draw (-1.5,0.866) node[anchor=south] {{\tiny $L^3$}};
\draw [stealth-stealth](0.5,-0.866) -- (-0.5,0.866);
\draw (0.5,-0.866) node[anchor=north] {{\tiny $S^{6}$}};
\draw (-0.5,0.866) node[anchor=south] {{\tiny $S^2$}};
\end{tikzpicture}
\end{figure} 

It follows from the formulas of $(c_1, c_2, c_3, c_4, c_5)$ given above that the vector fields $S^1$, $S^2$ and $S^3$ have the following coordinate expressions:
\begin{align}\label{s1}
S^1&=\frac{\ka c^2-h^2}{2c^3}\partial_\phi+\cos(\psi)(\frac{\ka c^2-h^2}{2c^3h}\partial_q-\frac{\partial_{\psi}}{h})-\sin(\psi) \partial_h,\\ \label{s2}
S^2&=\frac{\ka c^2-h^2}{2c^3}\partial_\theta+\sin(\psi)(\frac{\ka c^2-h^2}{2c^3h}\partial_q-\frac{\partial_{\psi}}{h})+\cos(\psi) \partial_h,\\ \label{s3}
S^3&=3(h^2-\ka c^2)(\theta \partial_{\phi}-\phi \partial_{\theta})-4c^2h(\ka \psi+2c q)\partial_h-(3 \ka c^2-h^2)\partial_\psi\\ \nonumber
&\quad{}-\frac{1}{2c^3}(c^4(2 c q+\ka \psi)^2-(3\ka c^2-h^2)(\ka c^2-h^2))\partial_q\\ \nonumber
&\quad{}+\cos(\psi)\bigg(6 c^3\phi \partial_h+\frac{3}{h}(h^2-c^2\ka)\theta \partial_q+\frac{6\theta c^3}{h}\partial_\psi-\frac{2h}{c}(\ka c^2-h^2)(2 cq+\ka \psi)\partial_{\theta}\\\nonumber
&\quad{}-\frac{h}{2c^3}(2c^4(2c q+\ka \psi)^2-(3\ka c^2-h^2)(\ka c^2-h^2))\partial_{\phi}\bigg)\\\nonumber
&\quad{}-\sin(\psi) \bigg(-6 c^3 \theta \partial_h+\frac{3}{h}(h^2-\ka c^2) \phi \partial_q+\frac{6c^3 \phi}{h}\partial_{\psi}-\frac{2h}{c}(\ka c^2-h^2)(2cq+\ka \psi)\partial_\phi\\\nonumber
&\quad{}+\frac{h}{2c^3}(2c^4(2 c q+\ka \psi)^2-(3\ka c^2-h^2)(\ka c^2-h^2))\partial_{\theta} \bigg).
\end{align}
The vector fields $S^1$ and $S^2$ are in the span of the distribution with $S^1=\frac{\ka c^2-h^2}{2c^3}X^2$ and $S^2=\frac{\ka c^2-h^2}{2c^3}X^1$ and they are annihilated by the 1-forms $\{\Theta_1, \Theta_2, \Theta_3\}$.

\begin{theorem}\label{s123}
The vector fields $S^1$, $S^2$ and $S^3$ given by the coordinate expressions (\ref{s1}), (\ref{s2}), (\ref{s3}) associated to the $(2,3,5)$-distribution of the configuration space of an An-Nurowski surface rolling without slipping or twisting over the plane, pairwise bracket generate the Lie algebra of split $\frak{g}_2$.
\end{theorem}    

The Lie algebra $\frak{g}_2$ is isomorphic to the symmetry algebra $\mathfrak{aut}(\cD)$, the Lie algebra associated to the group of automorphisms of $\cD$. The relationship with the coordinate expressions for the symmetry algebra can be found in \cite{r21c}. 
 The commutation relations for this Lie algebra can be found in the appendix to \cite{r21c}. The rest of the vector fields actually are not too bad to write down by computing them from the Lie brackets using the coordinate functions $(c_1,c_2,c_3,c_4,c_5)$. For example, we find that the commuting vector fields in the Cartan subalgebra are
 \begin{align*}
 H_1&=-h\partial_{h}-6\phi\partial_{\phi}-\frac{2cq+\ka \psi}{c}\partial_q+\frac{h}{2c^3}(3h^2-4\ka c^2)\cos(\psi) \partial_{\theta}+\frac{h}{2c^3}(3h^2-8\ka c^2)\sin(\psi)\partial_{\phi}\\
 &\quad{}+\sin(2\psi)(3\partial_{\psi}-\frac{3\ka}{2c}\partial_q)-3\cos(2\psi)h\partial_h-\frac{h^3}{2c^3}\sin(3\psi)\partial_{\phi}+\frac{h^3}{2c^3}\cos(3\psi)\partial_{\theta},\\
 H_2&=-6\theta \partial_{\theta}-6\phi \partial_{\phi}-2h \partial_{h}-\frac{2}{c}(2c q+\ka \psi) \partial_q+\frac{2h \ka}{c}\cos(\psi)\partial_{\theta}-\frac{2h \ka}{c}\sin(\psi) \partial_{\phi}.
 \end{align*}

We can also give the full list of vector fields in the case where $c=1$ and $\kappa=0$ (modulo homotheties, this corresponds to the class $g_{1o}$ in Theorem 3 of \cite{AN14}) as they can be displayed more concisely.
In the coordinates $(\theta, \phi, h,q,\psi)$, we have the short roots given by
\begin{align*}
S^1&=-\frac{h^2}{2}\partial_\phi-\cos(\psi)(\frac{h}{2}\partial_q+\frac{\partial_{\psi}}{h})-\sin(\psi) \partial_h,\\
S^2&=-\frac{h^2}{2}\partial_\theta-\sin(\psi)(\frac{h}{2}\partial_q+\frac{\partial_{\psi}}{h})+\cos(\psi) \partial_h,\\
S^3&=3h^2(\theta \partial_{\phi}-\phi \partial_{\theta})-8hq\partial_h+h^2\partial_\psi-\frac{1}{2}(4q^2-h^4)\partial_q\\
&\quad{}+\cos(\psi)\bigg(6 \phi \partial_h+3h\theta \partial_q+\frac{6\theta}{h}\partial_\psi+4h^3q\partial_{\theta}-\frac{h}{2}(8q^2-h^4)\partial_{\phi}\bigg)\\
&\quad{}-\sin(\psi) \bigg(-6 \theta \partial_h+3h \phi \partial_q+\frac{6\phi}{h}\partial_{\psi}+4h^3q\partial_\phi+\frac{h}{2}(8q^2-h^4)\partial_{\theta} \bigg),\\
S^4&=\partial_q+\sin(\psi) h \partial_{\theta}+\cos(\psi) h \partial_\phi,\\
S^5&=-(4 q^2-h^4)\partial_{\phi}+6\theta \partial_q+2qh^2\partial_{\theta}+\sin(\psi)(2h^2\partial_h+6\theta h \partial_{\theta}+\frac{8q}{h}(\partial_\psi-\frac{h^2}{4}\partial_q))\\
&\quad{}-\cos(\psi)(8q \partial_h-6 \theta h \partial_{\phi}-h^3\partial_q)+2\cos(2\psi)q h^2\partial_{\theta}-2\sin(2\psi) q h^2\partial_{\phi},\\
S^6&=-(4q^2-h^4)\partial_{\theta}-6\phi \partial_q-2qh^2\partial_{\phi}-\cos(\psi)(2h^2\partial_h+6\phi h \partial_{\phi}+\frac{8q}{h}(\partial_\psi-\frac{h^2}{4}\partial_q))\\
&\quad{}-\sin(\psi)(8q \partial_h+6 \phi h \partial_{\theta}-h^3\partial_q)+2\sin(2\psi) qh^2\partial_{\theta}+2\cos(2\psi)qh^2\partial_{\phi}.
\end{align*}
The long roots are given by the vector fields
 \begin{align*}
 L^1&=-\partial_{\theta},\\
 L^2&=\partial_{\phi},\\
 L^3&=3(2\theta \partial_{\phi}-\partial_\psi)-\frac{h^3}{2}(3\sin(\psi)+\sin(3\psi))\partial_{\theta}+\frac{h^3}{2}(3\cos(\psi)-\cos(3\psi))\partial_{\phi}\\
 &\quad{}+3\sin(2\psi)h\partial_h+3\cos(2\psi)\partial_\psi,\\
 L^4&=-((2q)(8q^2+3h^4)-36\phi \theta)\partial_{\phi}+(h^6+36\theta^2)\partial_{\theta}-18\phi \partial_{\psi}+18 \theta h \partial_h+36\theta  q\partial_q\\
 &\quad{}-3\sin(\psi)\bigg(8qh^2\partial_h+3h^3\phi\partial_{\theta}-3h^3\theta\partial_{\phi}-\frac{2}{h}(h^4+4q^2)\partial_{\psi}+8q^2h\partial_q\bigg)\\
 &\quad{}-\cos(\psi)\bigg(24q^2\partial_h-9h^3 \theta\partial_{\theta}-9 h^3 \phi\partial_{\phi}-6qh^3\partial_q\bigg)\\
 &\quad{}+3\sin(2\psi)(6\theta \partial_{\psi}+6h \phi \partial_{h}-4q^2h^2\partial_{\phi}+2qh^4\partial_{\theta})\\
 &\quad{}+3\cos(2\psi)(6 \phi \partial_{\psi}-6 h \theta \partial_{h}+4q^2h^2\partial_{\theta}+2qh^4\partial_{\phi})\\
 &\quad{}-3\cos(3\psi)h^3(\phi \partial_{\phi}-\theta \partial_{\theta})-3\sin(3\psi)h^3(\phi\partial_{\theta}+\theta \partial_{\phi}),\\
 L^5&=-((2q)(8q^2+3h^4)+36\phi \theta)\partial_{\theta}-(h^6+36\phi^2)\partial_{\phi}-18\theta \partial_{\psi}-18 \phi h \partial_h-36\phi q \partial_q\\
 &\quad{}+3\cos(\psi)\bigg(8h^2q\partial_h-3h^3\theta\partial_{\phi}+3h^3\phi\partial_{\theta}-\frac{2}{h}(h^4+4q^2)\partial_{\psi}+8q^2h\partial_q\bigg)\\
 &\quad{}-\sin(\psi)\bigg(24q^2\partial_h-9h^3 \theta\partial_{\theta}-9 h^3 \phi\partial_{\phi}-6qh^3\partial_q\bigg)\\
 &\quad{}-3\cos(2\psi)(6\theta \partial_{\psi}+6 h \phi \partial_{h}-4q^2h^2\partial_{\phi}+2qh^4\partial_{\theta})\\
 &\quad{}+3\sin(2\psi)(6\phi \partial_{\psi}-6 h \theta \partial_{h}+4q^2h^2\partial_{\theta}+2qh^4\partial_{\phi})\\
 &\quad{}+3\cos(3\psi)h^3(\phi \partial_{\theta}+\theta \partial_{\phi})+3\sin(3\psi)h^3(\theta \partial_{\theta}-\phi\partial_{\phi}),\\
 L^6&=3(2\phi \partial_{\theta}+\partial_\psi)-\frac{h^3}{2}(3\sin(\psi)+\sin(3\psi))\partial_{\theta}+\frac{h^3}{2}(3\cos(\psi)-\cos(3\psi))\partial_{\phi}\\
 &\quad{}+3\sin(2\psi)h\partial_h+3\cos(2\psi)\partial_\psi,
 \end{align*}
where in the vector fields we group the terms accordingly by the bases of trigonometric functions $(1,\cos(\psi), \sin(\psi), \cos(2\psi), \sin(2\psi), \cos(3\psi), \sin(3\psi))$. 
We also have
 \begin{align*}
 H_1&=-h\partial_{h}-6\phi\partial_{\phi}-2q\partial_q+\frac{3h^3}{2}\cos(\psi) \partial_{\theta}+\frac{3h^3}{2}\sin(\psi)\partial_{\phi}\\
 &\quad{}+3\sin(2\psi)\partial_{\psi}-3\cos(2\psi)h\partial_h-\frac{h^3}{2}\sin(3\psi)\partial_{\phi}+\frac{h^3}{2}\cos(3\psi)\partial_{\theta},\\
 H_2&=-6\theta \partial_{\theta}-2h \partial_{h}-4q\partial_q-6\phi \partial_{\phi}.
 \end{align*}

There is some symmetry in the coordinate expressions for the following pairs of vector fields: ($S^1$, $S^2$), ($S^5$, $S^6$), $(L^3,L^6)$ and $(L^4,L^5)$ on either side of the two hyperplanes separated by $S^3$ and $S^6$. Also note that the rotation term goes up to $2\psi$ in $S^5$ and $S^6$ while up to $3\psi$ in the long roots $L^3$, $L^4$, $L^5$, $L^6$.

The coordinates $(c_1,c_2,c_3,c_4,c_5)$ determine a natural projection of the Lie algebra of $\frak{sl}_3$ onto a 4-dimensional hypersurface in the following manner. Consider the subalgebra given by the $\R$-linear span of $\{L^1, L^2, L^3, L^4, L^5, L^6, H_1,H_2\}$. Referring to the coordinate expressions for the vector fields in the appendix of \cite{r21c}, only positive powers of $c_3$ appear in the coordinate expression and terms that contain $c_3\partial_{c_3}$. We can therefore restrict to the hypersurface $c_3=0$ to obtain the vector fields
\begin{align*}
L^1&=-6\partial_{c_1},\quad L^2=6\partial_{c_2},\quad L^3=6(c_1\partial_{c_2}-c_4\partial_{c_5}),\\
L^4&=6c_1(c_1\partial_{c_1}+c_2\partial_{c_2}+c_4\partial_{c_4})-6 c_2 c_4\partial_{c_5},\\
L^5&=-6c_2(c_1\partial_{c_1}+c_2\partial_{c_2}+c_5\partial_{c_5})+6 c_1 c_5\partial_{c_4},\\
L^6&=6(c_2\partial_{c_1}-c_5\partial_{c_4}),\\
H_1&=-6(c_2\partial_{c_2}-\frac{1}{3}c_4\partial_{c_4}+\frac{2}{3}c_5\partial_{c_5}),\\
H_2&=-6(c_1\partial_{c_1}+c_2\partial_{c_2}+\frac{1}{3}c_4\partial_{c_4}+\frac{1}{3}c_5\partial_{c_5})
\end{align*}
forming an $\frak{sl}_3$ algebra on $\{c_3=0\}$.

This parametrisation of $\frak{sl}_3$ comes from conformal rescalings of the split signature metric
\[
dS^2=2\der c_1 \der c_5+2 \der c_2\der c_4
\]
in $4$ dimensions. Let $c_6=c_1c_5+c_2c_4$.
We can rescale the metric by taking the conformal factor to be the inverse of one of the coordinates, say $c_4$, to obtain
\begin{align*}
\frac{1}{c_4^2}dS^2=2\der\left(\frac{c_6}{3 c_4}\right)\der\left(\frac{-3}{c_4}\right)+2\der\left(\frac{c_5}{c_4}\right)\der\left(\frac{c_1}{c_4}\right).
\end{align*}
This determines new coordinates given by
\begin{align*}
(c_1,c_2,c_4,c_5)\stackrel{\varphi}\mapsto \left(\frac{c_6}{3c_4},\frac{c_5}{c_4},\frac{c_1}{c_4},-\frac{3}{c_4}\right)=:(c_1^{(1)},c_2^{(1)},c_4^{(1)},c_5^{(1)}).
\end{align*}
The map $\varphi$ has the property that the 6-fold composition $\varphi^{(6)}$ is the identity map.
This corresponds to the cyclic group $C_6$ cycling through the six roots starting with $\partial_{c_1}$ and $\partial_{c_2}$ in the root diagram.
Together with the symmetry $(c_1,c_2,c_4,c_5)\mapsto(c_2, c_1,c_5,c_4)$, these two maps generate the dihedral group $D_6$ of order 12. 
The vector fields $\partial_{c_i}$, $\partial_{c_i^{(1)}}$, $\partial_{c_i^{(2)}}$, $\partial_{c_i^{(3)}}$, $\partial_{c_i^{(4)}}$, $\partial_{c_i^{(5)}}$ for $i=1$ or $2$, then determine up to constants $\{L^1, L^2, L^3, L^4, L^5, L^6\}$ above. The corresponding action of the symmetry group $SL_3(\R) \subset G_2$ can be compared with that acting on the dancing metric in Proposition 4.1 of \cite{BLN}, where $SL(3,\R)$ acts on $\R^3\times \R^3$, with the first factor projectivising to $\R\P^2$ and the second factor projectivising to its dual dancing pair $\R\P^{2*}$.
In the case of the circle twistor distribution of the An-Nurowski surface rolling on a plane, we have the change of coordinates on $\{c_3=0\}$ be given by
\begin{align*}
c_1&=6c^3\theta-h(3c^2\ka-h^2)\cos(\psi),\quad c_2=6c^3\phi+h(3c^2\ka-h^2)\sin(\psi),\\
c_4&=-\sin(\psi) h, \quad c_5=\cos(\psi) h. 
\end{align*}
What is more surprising, is that in this instance there is a further projection of $\frak{sl}_3$ onto the 3-dimensional submanifold $\R^2 \times \S^1$.

\section{Plane without a surface}

In the short story titled ``Solid Geometry" by the British novelist Ian McEwan, a fictional nineteenth century mathematician by the name of David Hunter presented a lecture about a plane without a surface in front of an initially skeptical mathematical audience.
In this section we consider the projection of the 5-dimensional circle twistor distribution associated to the An-Nurowski surface rolling on the plane to $\R^2 \times \S^1$. By restriction, we obtain an action of the symmetry group $SL(3,\R) \subset G_2$ on $\R^2 \times \S^1$, or the configuration space without the An-Nurowski surface as follows. 
\begin{theorem}
Let $\cD$ be the An-Nurowski $(2,3,5)$-distribution on the circle twistor bundle $M$ associated to the An-Nurowski surface $\Si$ with metric (\ref{anm1}) rolling without slipping or twisting on the plane. Let $\frak{g}_2\cong \frak{aut}(\cD)$ be the Lie algebra bracket-generated by the vector fields $S^1$, $S^2$ and $S^3$ in Theorem \ref{s123}. Then the vector fields of the subalgebra $\frak{sl}_3\subset \frak{g}_2$ naturally projects through 
\begin{align*}
\pi:M&\rightarrow \R^2 \times \S^1\\
(\theta,\phi,h,q,\psi)&\mapsto (\theta,\phi,0, -\frac{\ka}{2c}\psi,\psi)
\end{align*}
onto $\R^2\times \S^1=\{q=-\frac{\ka}{2c}\psi\}\cap\{h=0\}$, the submanifold of the configuration space without $\Si$. 
\end{theorem}
\begin{proof}
We consider the coordinate transformation of $(c_1,c_2,c_3,c_4,c_5)$ given in Section 3. We take the natural $\frak{sl}_3$ subalgebra acting on the 4-dimensional hypersurface $\{c_3=0\}=\{2c^3q+\ka c^2\psi=0\}$ as given above. We then observe now that in the long roots only positive powers of $h$ remain in the coordinate expressions and terms that contain $h\partial_h$. This allows us to project further the $\frak{sl}_3$ onto the submanifold  $\R^2\times \S^1=\{q=-\frac{\ka}{2c}\psi\}\cap\{h=0\}$. The $\S^1$ is the topology of the fibre coordinate of the circle twistor bundle. Taking for instance $c=1$, this gives the vector fields
\begin{align*}
L^1&=-\partial_{\theta},\\
L^2&=\partial_{\phi},\\
L^3&=-6\sin(\psi)^2\partial_{\psi}+6\theta \partial_{\phi},\\
L^4&=-36(\sin(\psi)\phi-\cos(\psi)\theta)\sin(\psi)\partial_{\psi}+36\theta(\theta\partial_{\theta}+\phi\partial_{\phi}),\\
L^5&=36(\sin(\psi)\phi-\cos(\psi)\theta)\cos(\psi)\partial_{\psi}-36\phi(\theta\partial_{\theta}+\phi\partial_{\phi}),\\
L^6&=6\cos(\psi)^2\partial_{\psi}+6\phi \partial_{\theta},\\
H_1&=6\cos(\psi)\sin(\psi)\partial_{\psi}-6\phi \partial_{\phi},\\
H_2&=-6\theta \partial_{\theta}-6\phi \partial_{\phi}
\end{align*}
 on $\R^2 \times \S^1$ which define the local action of $SL(3,\R)$. 
\end{proof}
We can contrast this result with the situation of Proposition 4.1 in \cite{BLN}, where $SL(3,\R)$ acts on the dancing pair. 
This 3-dimensional projection is peculiar to the example of An-Nurowski surface rolling on the plane and is not available for the scenarios of the maximally symmetric hyperboloids or spheres rolling distribution in \cite{r17} and \cite{r21b}.

Furthermore the metric representative $g$ as a symmetric 2-tensor vanishes identically on $\R^2 \times \S^1$, and it is mysterious what symmetries give rise to these $\frak{sl}_3$ generators on $\R^2\times \S^1$. 

We can also ask if there is a natural way to reverse the procedure and prolong the vector fields starting from the $\frak{sl}_3$ Lie algebra above to obtain the $\frak{g}_2$ Lie algebra. The expressions involving the prolongation coordinates $\partial_h$ and $\partial_q$ are not so straightforward, from the displayed formulas above. 

The three $\frak{sl}_2$'s in the $\frak{sl}_3$ algebra generated by the antipodal long roots can be written in a dual basis of 1-forms of rank $3$, and we can further use the construction of $SL(2)$ Pfaffian systems in \cite{r21a} to create more differential systems with split $\frak{g}_2$ symmetry. However, the geometric link with the circle twistor construction using the original An-Nurowski surface would then be lost.

 \end{document}